\documentclass{article}
\usepackage{graphicx} 
\usepackage[letterpaper, margin=1in]{geometry}
\usepackage{authblk} 
\usepackage{amssymb}
\usepackage{amsmath}
\title{\textsc\slshape A Complete Derivation of Complex Circle Manifold (CCM) Riemannian manifold Optimization Equations}
\author[I]{Amirreza Tabrizi}
\author[II]{Mohammad Hadi Mirmohammadi
}

\affil[I]{\small School of Electrical Engineering and Computer Engineering, University of Tehran

\small\textsc{\slshape amtabrizi@ieee.org}
}
\affil[II]{\small School of Electrical Engineering, Iran University of Science and Technology 

\small\textsc{\slshape mo\_mirmohammady1383@elec.iust.ac.ir}
}

\date{August 2025}

\begin{document}

\maketitle
\section*{\centering Abstract}
After reviewing manifold optimization techniques in applications like MIMO communication systems, phased array beamforming\footnote{for further reading see[1]}, radar, control theory, we observed that the Complex Circle Manifold (CCM) is widely employed, yet its foundational relations and equations lack a rigorous, self-contained derivation in the literature. This paper provides a systematic and rigorous proof of CCM’s key properties, including its tangent space and Riemannian gradient operations, with explicit connections to real-world optimization problems. Our work aims to serve as a unified reference for researchers and practitioners applying CCM Manifold Optimization.

\section*{\centering Introduction} 
Modern wireless communication systems and array processing techniques frequently employ complex-valued beamforming vectors constrained to lie on the complex circle manifold (CCM), where each element maintains unit modulus ($|w_i| = 1$). This constant modulus constraint enables uniform power distribution while facilitating efficient hardware implementation in phased array systems. However, optimizing over this constrained manifold presents significant theoretical and computational challenges that cannot be adequately addressed using conventional Euclidean optimization methods.

The inherent geometric structure of the CCM naturally calls for Riemannian optimization techniques, which properly account for the manifold's curvature through tangent space projections and retraction operations. These methods preserve the unit modulus constraint throughout the optimization process. A particularly challenging aspect arises from the need to optimize real-valued cost functions with respect to complex parameters. Optimizing real-valued functions over complex manifolds introduces fundamental difficulties due to the inherent properties of complex differentiation. Unlike real-valued calculus, complex differentiability requires analyticity (holomorphicity) via the Cauchy-Riemann equations, which most real-valued cost functions fail to satisfy. This non-analytic nature means standard complex derivatives do not exist, forcing alternative approaches like Wirtinger calculus or mapping to $\mathbb{R}^2$.
To simplify complex optimization, we employ ${\mathbb{C}}{\mathbb{R}}$-calculus\footnote{for further reading see[2]} by establishing a bijective mapping from the complex field $(\mathbb{C}, +, \cdot)$ to $(\mathbb{R}^2, +)$ vector space. This transformation removes the multiplicative structure while preserving the additive one, enabling conventional real-valued calculus. The mapping $z = x + jy \leftrightarrow (x,y)$ allows us to:
\begin{itemize}
\item Compute derivatives using real analysis
\item Maintain the complex domain's geometric interpretation
\item Derive optimization formulations consistent with literature
\end{itemize}
While existing literature frequently cites Book [3] for derivations of Complex Circle Manifold (CCM) equations, these references often lack complete mathematical rigor and detailed proofs. To address this gap, we present in the following sections a comprehensive and rigorous derivation of these fundamental equations, establishing each result with careful mathematical justification.

\section*{\centering{Basics}}
We now examine one of the most frequently formulated optimization problems involving the Complex Circle Manifold:

$$\min_{x \in \mathcal{M}}{{x^H}{A}{x}} $$
where the real-valued cost function $f(x) = x^H A x$ , $f :\mathcal{M}\subseteq\mathbb{C}^n\to \mathbb{R}$ is defined for complex vectors 
$x \in \mathbb{C}^n$ constrained to lie on the embedded manifold 
$\mathcal{M} \subseteq \mathbb{C}^n $, with $A \in \mathbb{C}^{n \times n}$ being 
a Hermitian matrix.
now lets briefly look at the $\mathbb{C}\mathbb{R}$-calculus we will use :
the field of complex numbers, $\mathbb{C}:\{ (x,y)\mid  x,y \in \mathbb{R} \}$ contains multiplicative structure $\cdot_{\mathbb{C}} : \mathbb{C} \times \mathbb{C} \to \mathbb{C}$  and additive structure $+_{\mathbb{C}} : \mathbb{C} \times \mathbb{C} \to \mathbb{C}$.
now if we restrict the domain of multiplication to $D=(R \times \mathbb{C}) \cup (\mathbb{C} \times R)$ in which $R=\{(x,y) \in \mathbb{C} \mid y=0 \}$.obviously it has an algebra-preserving bijective map with  vector space $\mathbb{R}^2$ over field $\mathbb{R}$ ( which resembles to an isomorphism)with standard additive and multiplicative structures.
$$(\mathbb{C}, +_\mathbb{C}, \cdot_\mathbb{C} |_D)  \cong (\mathbb{R}^2 , +,\cdot)$$
thus from now on we will do our problem on  equivalent $\mathbb{R}^{2n}$ .
now by defining a standard inner product ${<}{\cdot,\cdot}{>}$ , we can work in riemannian geometry .
suppose
$v =\sum_{i} {v_i}{e_i} $
 and
$w =\sum_{i} {w_i}{e_i} $ 
. then the inner product will be defined as :
$${<}{v,w}{>} \overset{\Delta}{=}\sum_{i} {v_i}{w_i} $$

\section*{\centering Derivative}
first we will expand $x^{H}Ax$ . suppose $A=(a_{dk})_{n \times n}$  , $x=(w_{i})_{n \times 1}$ and $w_i= {x_{i}}^{1} + j{x_{i}}^{2}$ ( j is the imaginary unit). after doing needed multiplications we will have :
$$x^H A x = \sum_{i=1}^n \sum_{k=1}^n \overline{w}_i a_{ik} w_k$$
 since A is hermitian the imaginary part will cancel out and the remaining after more expansion is :
 $$x^H A x = \sum_{i=1}^n \sum_{k=1}^n \left( x_i^1-jx_i^2 \right) a_{ik} \left( x_k^1+jx_k^2 \right)$$
 now lets take its partial derivative respect to each ${x_m^l}$ , $l=1,2$ . 
 $$ \frac{\partial }{\partial {x_m^1} } x^{H}Ax =\frac{\partial }{\partial {x_m^1} } \sum_{i=1}^n \sum_{k=1}^n \left( x_i^1-jx_i^2 \right)\left( x_k^1+jx_k^2 \right) a_{ik}$$
 $$ = \sum_{i=1}^n \sum_{k=1}^n  \left(\frac{\partial }{\partial {x_m^1} }\left( x_i^1-jx_i^2 \right)\left( x_k^1+jx_k^2 \right) + \left( x_i^1-jx_i^2 \right)\frac{\partial }{\partial {x_m^1} }\left( x_k^1+jx_k^2 \right)\right) a_{ik}$$
 $$\sum_{i=1}^n \sum_{k=1}^n  \left(\delta_{mi}w_k + \overline{w}_i \delta_{mk} \right)a_{ik} = \sum_{i=1}^n \sum_{k=1}^n  \delta_{mi}w_k a_{ik} + \sum_{i=1}^n \sum_{k=1}^n  \overline{w}_i \delta_{mk}a_{ik}$$
 $$\sum_{k=1}^n  w_k a_{mk} + \sum_{i=1}^n  \overline{w}_i a_{im} = \sum_{k=1}^n  w_k a_{mk} + \overline{\left(\sum_{i=1}^n  w_i \overline{a}_{im}\right)}$$
  $$\sum_{k=1}^n  w_k a_{mk} + \overline{\left(\sum_{i=1}^n  w_i a_{mi}\right)} = 2\operatorname{Re}\left\{\sum_{k=1}^n  w_k a_{mk}\right\}$$
and partial derivative respect to $x_m^2$
$$\frac{\partial }{\partial {x_m^2} } x^{H}Ax = \sum_{i=1}^n \sum_{k=1}^n  \left(-j\delta_{mi}w_k  +j \overline{w}_i \delta_{mk} \right)a_{ik}$$
$$-j\left( \sum_{k=1}^n  w_k a_{mk} - \overline{\left(\sum_{i=1}^n  w_i \overline{a}_{im}\right)}\right)=-j\left(2j\operatorname{Im}\left\{\sum_{k=1}^n  w_k a_{mk}\right\}\right)$$
$$2\operatorname{Im}\left\{\sum_{k=1}^n  w_k a_{mk}\right\}$$

now we define complex partial derivatives respect to $w_m$ :
$$\frac{\partial}{\partial w_m} \overset{\Delta}{=} \frac{\partial}{\partial x_m^1} + j\frac{\partial}{\partial x_m^2}   $$
 using above definition , now we are able to define the complex gradient operator :
 $$\nabla_w \overset{\Delta}{=}\left(\frac{\partial}{\partial w_m}\right)_{n\times1}    $$

now if we take the gradient of $F=x^H A x$ respect to $w$
we will have :

$$\nabla_{w}f = \begin{pmatrix}\frac{\partial f}{\partial w_1} & \frac{\partial f}{\partial w_2}  &   \hdots & \frac{\partial f}{\partial w_n}\end{pmatrix}^T =2 \begin{pmatrix}\sum_m{a_{1m}}{w_m} & \sum_m{a_{2m}}{w_m}& \hdots & \sum_{m}{a_{nm}}{w_m} \end{pmatrix}^T$$
$$ = 2{A}{x}$$
\section*{\centering Complex Circle Manifold}

complex circle manifold is a product manifold $\mathcal{M}$ of $\Omega=\{ z \in \mathbb{C} \mid {|z|=1}\}$ . so $\mathcal{M}=\prod_{i=1}^{n} \Omega$.
\section*{Tangent space}
 for two points on product manifold $\mathcal{B} \times \mathcal{N}$ embedded in $\mathbb{R}^n \times \mathbb{R}^m$ like p and q we know the following relation between tangent spaces \footnote{for further reading and proof see[4]} :
$$T_{p \times q}({\mathcal B} \times {\mathcal {N}}) \cong {T_{p}{\mathcal{B}}} \oplus {T_{q}{\mathcal{N}}}$$
$$and$$
$$(T_{p \times q}({\mathcal B} \times {\mathcal {N}}))^{\perp} \cong ({T_{p}{\mathcal{B}}})^{\perp} \oplus ({T_{q}{\mathcal{N}}})^{\perp}$$
thus :
 $$ dim(T_{p \times q}({\mathcal B} \times {\mathcal {N}}))=dim({T_{p}{\mathcal{B}}})+dim({T_{q}{\mathcal{N}}})$$

and according to Rank-theorem :
 $$dim(T_{p \times q}({\mathcal B} \times {\mathcal {N}}))+dim((T_{p \times q}({\mathcal B} \times {\mathcal {N}}))^\perp)=n+m $$

for point x in the product manifold $\mathcal{M}$ :
 $$dim(T_{x}{\mathcal{M}})= \sum_{i=1}^{n} dim({T_{w_i}{\Omega}})=n \times 1 = n$$

thus $dim((T_{x}{\mathcal{M}})^\perp)=n$. 
we also know :
$$(T_x{\mathcal {M}})^\perp \cong \bigoplus_{i=1}^n (T_{w_i}{\Omega})^\perp$$
now our purpose is to determining basis for perpendicular space of manifold. to achieve this objective, in the first step ,we will select one non-trivial unit vector from $(T_{w_1}{\Omega})^\perp$ and n-1 trivial zero vector from $((T_{w_i}{\Omega})^\perp)_{i=2}^n$ and in the second step we will repeat this process by choosing that unit vector from $(T_{w_2}{\Omega})^\perp$ and choosing other trivial zero vectors from the other tangent spaces . by repeating this process n times , we actually constructed a basis for $(T_x{\mathcal{M}})^\perp$.
so we will say an arbitrary vector $z \in T_x{\mathcal{M}}$ if for all $b_i$ as the i-th basis for the $(T_x{M})^\perp$  $<z , b_i>=0$ .  after accomplishing this goal we should retrieve our $b_i$ in complex form $O_i$ so we could do our desired computations .

Let us mathematically find the perpendicular vector of $\Omega$. To do this an algebraic  equation is needed. using the levels set of our manifold  $|w_i| = 1$ we will have :
$$ |w_i|=\sqrt{(x_i^1)^2 + (x_i^2)^2} = 1$$
so the algebraic  equation that characterize manifold is : 
$$ \varphi=(x_i^1)^2 + (x_i^2)^2 -1 =0$$
we have to find the gradient of the zero level set :
$$\nabla \varphi= \begin{pmatrix}\frac{\partial \varphi}{\partial x_i^1} & \frac{\partial \varphi}{\partial x_i^2}\end{pmatrix}^T=\begin{pmatrix}2x_i^1  & 2x_i^2\end{pmatrix}^T$$
we will define $n_i$ as unit perpendicular vector of $\Omega$ as :
$$n_i \overset{\Delta}{=} \frac{\nabla \varphi}{|\nabla \varphi|} = \begin{pmatrix}x_i^1  & x_i^2\end{pmatrix}^T$$
the complex form of the above unit perpendicular vector is :
$$d_i \overset{\Delta}{=} \frac{\frac{\partial \varphi}{\partial w_i}}{|\frac{\partial \varphi}{\partial w_i}|}=\frac{\frac{\partial \varphi}{\partial x_i^1}}{|\frac{\partial \varphi}{\partial x_i^1}|} + j \frac{\frac{\partial \varphi}{\partial x_i^2}}{|\frac{\partial \varphi}{\partial x_i^2}|}$$
 now we should retrieve it in our desired complex form $O_i$ :

$$ O_i {=} \begin{pmatrix} 0 \\ 0 \\ \vdots \\d_i \\ \vdots \\ 0  \end{pmatrix} _{n \times 1}{=}\begin{pmatrix} 0 \\ 0 \\ \vdots \\w_i \\ \vdots \\ 0  \end{pmatrix}_{n \times 1}$$

note that we have two representations of vectors ,the euclidean one :

$$v_{euclidean} = \begin{pmatrix} v_1^1 & v_1^2 & v_2^1 & v_2^2 & \hdots & v_n^1 & v_n^2 \end{pmatrix}_{1 \times 2n} ^T$$
 and the complex form :
 $$v_{complex} = \begin{pmatrix} v_1^1 + jv_1^2 & v_2^1 + jv_2^2 & \hdots & v_n^1 + jv_n^2 \end{pmatrix}_{1\times n} ^T$$

the equivalence between them is result of excluding multiplicative structure from the complex algebra and because of this every $v_i^1$ and $v_i^2$ can reciprocally act on its alias in other vector form (in the case of multiplicative structure's existence , the equivalence  form was a matrix not a vector) .

we say vector $z \in T_x{\mathcal{M}}$ if for all basis $b_i\in (T_x{\mathcal{M}})^\perp$ the inner product 

$$<z,b_i>={z_i^1}{n_i^1} + {z_i^2}{n_i^2}={z_i^1}{x_i^1} + {z_i^2}{x_i^2}=0$$

consider $\operatorname{Re}\{z \odot \overline x \}= \operatorname{Re}\{\begin{pmatrix} (z_i^1+jz_i^2)(x_i^1-jx_i^2)  \end{pmatrix} _{n \times 1} \}$

$$= \operatorname{Re}\{\begin{pmatrix} {z_i^1}{x_i^1} +{z_i^2}{x_i^2} +j({z_i^2}{x_i^1}-{z_i^1}{x_i^2}) )  \end{pmatrix} _{n \times 1}\}=
\begin{pmatrix} {z_i^1}{x_i^1} +{z_i^2}{x_i^2}  \end{pmatrix}_{n \times 1} = 0_n $$ 
in which $\odot$ is \textbf{Hadamard} element-wise matrix product. and this is the same condition for perpendicularity of vector .

finally definition of our manifolds tangent space is :

$$T_x{\mathcal{M}}=\{z \in \mathbb{C}^n \mid \operatorname{Re}\{z \odot \overline x \}=0_n  \}$$

\section*{projection}
projection $\mathcal{P}_x{(z)}$ of vector z , on tangent space of point x on manifold $\mathcal{M}$ is :
$$\mathcal{P}_x{(z)}\overset{\Delta}{=}z-\sum_{b_i \in T_x{\mathcal{M}}}\frac {<z,b_i>}{<b_i,b_i>}b_i =z-\sum_{b_i \in T_x{\mathcal{M}}}{<z,b_i>b_i }$$
$$=z - \sum_{b_i \in T_x{\mathcal{M}}}({{z_i^1}{x_i^1} + {z_i^2}{x_i^2}})b_i=z- \sum_{b_i \in T_x{\mathcal{M}}}\begin{pmatrix} 0 \\ \vdots\\ z_i^1 (x_i^1)^2 + z_i^2x_i^1x_i^2 \\z_i^2 (x_i^2)^2 + z_i^1x_i^1x_i^2 \\0\\ \vdots \\0\end{pmatrix}_{2n \times 1}$$

$$=z-\begin{pmatrix} z_1^1 (x_1^1)^2 + z_1^2x_1^1x_1^2 \\ z_1^2 (x_1^2)^2 + z_1^1x_1^1x_1^2\\\vdots\\ z_i^1 (x_i^1)^2 + z_i^2x_i^1x_i^2 \\z_i^2 (x_i^2)^2 + z_i^1x_i^1x_i^2 \\ \vdots \\z_n^1 (x_n^1)^2 + z_n^2x_n^1x_n^2\\z_n^2 (x_n^2)^2 + z_n^1x_n^1x_n^2\end{pmatrix}_{2n \times 1}$$

now consider $\operatorname{Re}(z\odot \overline x)\odot x=\begin{pmatrix} {z_i^1}{x_i^1} +{z_i^2}{x_i^2}  \end{pmatrix} _{n \times 1} \odot \begin{pmatrix} x_i^1 + jx_i^2 \end{pmatrix}_{n \times 1}$

$$= \begin{pmatrix}  z_i^1 (x_i^1)^2 + z_i^2x_i^1x_i^2 + j(z_i^2 (x_i^2)^2 + z_i^1x_i^1x_i^2)  \end{pmatrix}_{n \times 1} $$

as we can see this is the complex representation of the former calculation . 
so we can restate the projection like this :

$$\mathcal{P}_x{(z)}=z-\operatorname{Re}\{z\odot \overline x\}\odot x$$

finally we can define the most important operator in manifold optimization , the  riemannian gradient :
$$\nabla_{\mathcal{M}}f \overset{\Delta}{=} \mathcal{P}_x({\nabla{f}})=\nabla_{w}f-\operatorname{Re}\{\nabla_{w}f\odot \overline x\}\odot x$$

\renewcommand{\refname}{}

\end{document}